\documentclass[12pt]{article}

\usepackage{latexsym}
\usepackage{amsfonts}
\usepackage{amsmath}
\usepackage{amssymb}

\newcommand{\op}[1]{\operatorname{#1}}

\renewcommand{\O}{{\mathcal O}}
\renewcommand{\P}{{\mathbb P}}

\newcommand\A{{\mathbb A}}
\newcommand\Z{{\mathbb Z}}
\newcommand\Q{{\mathbb Q}}
\newcommand\C{{\mathbb C}}

\newcommand\G{{\mathbb G}}
\newcommand\xh{{\hat{x}}}
\newcommand\yh{{\hat{y}}}

\newcommand\spec{{\mathsf{Spec}\,}}

\newcommand\D{{\mathcal D}}
\renewcommand\spec{{\mathsf{Spec}\,}}

\renewcommand\k{{\bf k}}

\newtheorem{dfn}{Definition}

\newtheorem{thm}{Theorem}
\newtheorem{prop}{Proposition}
\newtheorem{conj}{ Conjecture}

\title{Holonomic $\D$-modules and positive characteristic}

\begin{document}

\author{Maxim Kontsevich}
\maketitle

\section{Quantization and $\D$-modules}

\subsection{Quantization}

Let us recall the heuristic dictionary of quantization.
 $$\begin{array}{rcl}
\mbox{ symplectic $C^{\infty}$-manifold }(Y,\omega)& \longleftrightarrow & \mbox{ complex Hilbert space } H\\
\mbox{ complex-valued functions on } Y& \longleftrightarrow & \mbox{  operators in } H\\
\mbox{Lagrangian submanifolds } L\subset Y& \longleftrightarrow & \mbox{ vectors } v\in H\\
\end{array}$$

In the basic example of the cotangent bundle $Y=T^*X$ the space $H$ is  $L^2(X)$, functions on $Y$ which are polynomial
 along fibers of the projection $Y=T^*X\to X$ correspond to differential operators, and with the Lagrangian manifold $L\subset Y$
 of the form $L=\op{graph} dF$ for some function $F\in C^\infty(Y)$ we associate  (approximately) vector
  $\exp(iF/\hbar)$ where $\hbar\to 0$ is a small parameter (``Planck constant'').

Our goal in these lectures is to give an evidence for a hypothetical analog of the quantization in the algebraic case, based on the 
reduction to positive characteristic, briefly mentioned at the end of 
\cite{KaKo}. The idea to study quantization in positive characteristic was also used in \cite{Tsuchimoto}, and in the fundamental article 
\cite{OgusVologodsky}. It turns out that the correspondence 
between the classical and the quantum in the algebraic case is not one-to-one, but only between certain natural
 families of Lagrangian manifolds (or cycles) and of holonomic $\D$-modules, closely related to integrable systems. We formulate a series
 of conjectures about this correspondence.
  In the case of one variable (which in a sense contain all the keys to the general case) one can make all the constructions explicit and elementary.

\subsection{Reminder about holonomic $\D$-modules}

In the algebraic setting there is no obvious analog of the Hilbert space,
 even in the case of the cotangent bundle.
 A possible replacement for the notion of a function is the one of a holonomic $\D$-module.
Here we recall the definition and  several basic and well-known  facts (the standard reference is \cite{Borel}).

Let $X$ be a smooth affine algebraic variety over field $\k$ of zero characteristic, $\dim X=n$.
The ring $\D(X)$ of differential operators is $\k$-algebra of operators acting on $\O(X)$, generated by
functions and derivations:
$$f\mapsto gf,  \,\,f\mapsto \xi(f),\,\,\,g\in \O(X),\,\xi \in \Gamma(X,T_{X/\spec \k})\,\,.$$
Algebra $\D(X)$ carries the
 filtration $\D(X)=\cup_{k\ge 0}\D_{\le k}(X) $ by the degree of operators,
the associated graded algebra is canonically isomorphic to the algebra of functions on $T^*X$. In geometric terms, the grading comes
from the dilation by $\G_m$ 
 along the fibers of the cotangent bundle.

Let $M$ be a finitely generated $\D(X)$-module, and choose a finite-dimensional subspace $V\subset M$ generating $M$.
 Then consider the filtration 
$$M_{\le k}:=\D_{\le k}(X)\cdot V\subset M,\,\,\,k\ge 0\,\,. $$
The associated graded module $\op{gr}(M)$ is a finitely generated $\O(T^*X)$-module.
 By theorem of Gabber, its support (which is  a reduced conical subscheme in $T^*X$)
$$\mathsf{supp}\, (\op{gr}(M))\subset T^*X$$
 is coisotropic. In particular, the dimension of any irreducible
component  is at least $n=\op{dim} X$. The support does not depend on the choice of the generating space $V$, and is denoted by $\mathsf{supp}\,(M)$. 

A finitely generated module $M$ is called {\it holonomic} if and only if the dimension of its support is exactly $n$.
Let $(L_i)_{i\in I}$ be the set of irreducible components of $\mathsf{supp}\, (M)$. Each $L_i$ is the generic point of a conical Lagrangian
subvariety of $T^*X$. One can show that the multiplicity $m_i\ge 1$ of  $\op{gr}(M)$ at  $L_i$ does not depend
on the choice of generators as well. The notion of support with multiplicities is covariant with respect to
 automorphisms of $X$ (or, more generally, contravariant with respect to  \'etale maps).

 For general smooth $X$, not necessarily affine,
 we have a sheaf (in Zariski topology)
$\D_X$ of algebras
 of differential operators. Then 
 one can define $\D_X$-modules (and also holonomic modules) and their support. Finitely generated $\D_X$-modules form a  noetherian abelian category. For any two holonomic $\D_X$-modules $M_1,M_2$ we have
 $$
 \begin{matrix}
 \dim \op{Ext}^i(M_1,M_2)<\infty&\forall i \in \Z, \hfill \cr
 \op{Ext}^i(M_1,M_2)=0 \hfill &\mbox { for }i<0 \mbox{ or } i>n=\dim X\,.
 \end{matrix}
 $$
 Any algebraic vector bundle (a locally trivial sheaf of $\O_X$-modules)
 ${\cal E}\in Coh(X)$ endowed with a flat connection $\nabla$ has a natural structure of a holonomic $\D_X$-module. Its support (with multiplicities) is  the zero section of $T^*X$ taken with the multiplicity equal to $\op{rank}({\cal E})$.
 For any holonomic $\D_X$-module $M$ there exists a non-empty Zariski open subset $U\subset X$
  such that the restriction of $M$ is a bundle with a flat connection.

In the case of the affine space $X=\A^n_\k$ one can use another filtration called the Bernstein filtration. Namely, the algebra $\D(\A^n_\k)$ is the $n$-th
Weyl algebra $A_{n,\k}$
over\footnote{Same formulas  give the definition of  algebra $A_{n,R}$  for arbitrary commutative ring $R$.} $\k$, i.e., it has the presentation
  $$\k\langle \xh_1,\dots, \xh_{2n}\rangle/([\xh_i,\xh_j]=\omega_{ij},\,\,\,\forall\,i,j\,\,\, 1\le i,j\le 2n)\,,$$
where 
$$\omega_{ij}=\delta_{i,n+j}-\delta_{j,n+i}\,\,.$$
The generators are realized as 
$$\xh_i=x_i,\,\,\xh_{n+i}=\partial/\partial x_i,\,\,1\le i\le n\,.$$ 
The filtration is given by the degree in generators.
In this case the notion of a holonomic module is the same, but the support is now a Lagrangian cone in $\A^{2n}$,
 invariant with respect to the total dilation in $\A^{2n}$. Notice that the notion of a holonomic module over $A_{n,\k}$  and its support
is covariant with respect to the action of the symplectic group $Sp(2n,\k)$.

\subsection{Lagrangian cycles}

Let $(Y,\omega)$ be
a symplectic smooth quasi-projective algebraic variety  over a field $\k$ of arbitrary characteristic.
 \begin{dfn} A {\it Lagrangian subvariety}
  in $Y$ is the generic point of a smooth (not necessarily closed) irreducible Lagrangian submanifold $L\subset Y$. An {\it effective Lagrangian cycle} is a formal linear combination with coefficients in $\Z_{\ge 0}$
  of
  Lagrangian subvarieties.
 \end{dfn}
 Effective Lagrangian cycles form a submonoid in the abelian group of algebraic cycles
  $Z^n (Y)$ where $n:=\dim Y/2$.
 
 Let us assume that $Y$ is quasi-projective, fix a projective embedding $i:Y \hookrightarrow \P^N_\k$ and a constant $d>0$. Then
the set of effective Lagrangian cycles of degree $\le d$ with respect to $i$
   has a natural structure of the set of $\k$-points
   of a constructible set defined over $\k$. If we do not bound the degree, we obtain
    an ind-constructible set over $\k$ which we will denote by $ELC(Y)$.
In the case  $\mathsf{char}(\k)=0$ it is not 
inconceivable that $ELC(Y)$ is in fact an ind-scheme, not just merely an ind-constructible set.
Also, the same could be true for Lagrangian cycles of a bounded degree when $\mathsf{char}(\k)$ is large enough.
In section 3.1 we present some evidence for it.
    
    In the case $\mathsf{char}(\k)=0$ the support (with multiplicities) of a holonomic $\D_X$-module is
     a $\k$-point of $ELC(T^*X)$ invariant under the dilation by $\G_m$ along fibers  
     of  the bundle $T^*X\to X$. In the case $X=\A^n_\k$ and of the Bernstein filtration the support
      is invariant under the total dilation.

\section{Reduction of $\D$-modules to positive characteristic}

\subsection{Differential operators in positive characteristic}

Let $X$ be a smooth affine scheme over an arbitrary commutative ring $R$.
 We define the ring $\D(X)$ of
differential operators on $X$ as an $R$-linear associative algebra  generated additively by  $\O(X)$ and  by $T(X):=\Gamma(X,T_{X/\spec  R})$ subject to the following relations:
$$f\cdot f'=ff',\,\,\,f\cdot \xi= f\xi, \,\,\,\xi\cdot f-f\cdot \xi=\xi(f),\,\,\,\xi\cdot \xi'-\xi'\cdot\xi=[\xi,\xi'] $$
where $f,f'\in \O(X),\,\,\xi,\xi'\in T(X)$.

This definition is equivalent to the usual one if $R=\k$ is a field of characteristic zero. 
In general,  $\D(X)$ maps to the algebra of differential operators in the sense of Grothendieck, but
 the map is neither surjective, nor injective.
The surjectivity fails e.g. for 
$$R=\Z,\,\,X=\A^1_R,\,\,\,\O(X)=\Z[x]\,.$$
In this  case  the divided powers of $d/dx$ acting on $\O(X)$ 
$$\frac{(d/dx)^n}{n!}: x^m\mapsto \binom{m}{n} x^{m-n}$$
do not belong to the image of $\D(X)$ for $n\ge 2$.

The injectivity fails in positive characteristic:
$$R=\Z/p \Z,\,\,\, Y=\A^1_R, \,\,\,(d/dx)^p\mapsto 0\in \op{End}(\O(X))\,.$$

In general, for $R$ in characteristic $p>0$, i.e., $p\cdot 1_R=0$, the  algebra $\D(X)$ has a big center: 
$$\op{Center}(\D(X))\simeq \O(T^*X')$$
where $X'/\spec  R$ is the pullback of $X$ under the Frobenius map
 $$\op{Fr}_{R}:\spec R\to \spec  R,\,\,\,\op{Fr}^*_{R}(a)=a^p\,.$$
Moreover, the algebra $\D(X)$ is an Azumaya algebra of its center, it is a twisted form of the matrix algebra
  $\op{Mat}(p^n\times p^n,\O(T^* X'))$ where $n=\dim X$. In the basic example of the Weyl algebra
  $A_{n,R}:=\D(\A^n_R)$, the center is the algebra of polynomials
  $R[\xh_1^p,\dots,\xh_{2n}^p]$.

\subsection{Infinitely large prime}

 It will be convenient to introduce the following notation
(``reduction modulo infinitely large prime'') for an arbitrary commutative ring
  $R$:
$$
  R_\infty:=\varinjlim_{f.g.\,\widetilde{R}\subset
  R}  \left(\,\,\prod_{\mathrm{primes}\,\,p} (\widetilde{R}\otimes \Z/p \Z)
  \,\,\,\,\big{/} \, \bigoplus_{\mathrm{primes}\,\,p} (\widetilde{R}\otimes \Z/
  p\,\Z)\right)\,\,.
  $$

  Here the inductive limit is taken over the filtered system consisting
  of all finitely generated subrings $\widetilde{R}
\subset R$, and the index $p$ runs over
  primes $2,3,5,\dots$.
It is easy to see that the ring $R_\infty$ is defined over $\Q$ (all
primes are invertible in $R_\infty$), and the obvious map $R\mapsto
R_\infty$ gives an inclusion $R\otimes \Q\hookrightarrow R_\infty$.
Also, there is a universal Frobenius endomorphism given by $$\op{Fr}_{R_\infty}^*:R_\infty\to
R_\infty,\,\,\,\op{Fr}_{R_\infty}^*(a_p)_{\,\mathrm{primes}\,\,p}:=(a_p^p)_{\,\mathrm{primes}\,\,p}\,\,.$$

Now we discuss a related notion.
Let $R=\k$ be a field of characteristic zero, and $S/\k$ be a constructible set.
We define the set of ``twisted points modulo large primes'' as
 $$ S_\infty^{tw}:= \varinjlim_{f.g.\,\widetilde{R}\subset
  \k, \,\widetilde{S}} \left(\,\,\prod_{\mathrm{primes}\,\,p} \left(C.S.:
\widetilde{S}'_{\widetilde{R},p}\to
\spec\widetilde{R}\otimes \Z/p \Z\right)
  \right)\,\,, $$
where the limit is taken over pairs consisting of a finitely generated subring $\widetilde{R}\subset \k$ and an affine scheme of finite type
 $\widetilde{S}/\widetilde{R}$
 endowed with a constructible  over $\k$ bijection between $\widetilde{S}\times_{\spec \widetilde{R}}\spec\k$ and $S$.
The scheme $\widetilde{S}'_{\widetilde{R},p}$  is defined as 
 the pullback under $\op{Fr}_{\spec(\widetilde{R}\otimes \Z/p \Z),p}$ of the scheme $\widetilde{S}\times \spec\Z/p\,\Z$.
Finally, the abbreviation $C.S.$ means the set of constructible sections of a morphism of schemes of finite type.

 Notice that we automatically identify collections of constructible sections which differ only at a finite set of primes. The reason is that $R=\k$ contains $\Q$ and we can always add to $\widetilde{R}$ inverses of any finite set of primes.

In the special case $S=\A^1_\k$ we have an embedding
$$\k_\infty\hookrightarrow (\A^1_\k)_\infty^{tw}\,.$$
For an ind-constructible set over $\k$ represented as a countable limit of constructible sets $S=\displaystyle{\lim_{\longrightarrow}}\, S_i$
 we define $ S_\infty^{tw}$ as the inductive limit of sets $ (S_i)_\infty^{tw}$.

\subsection{Reduction of finitely generated $\D$-modules}

Let $X$ be a smooth affine\footnote{All the considerations here extend to the case of not necessarily affine $X$.} 
variety over field $\k, \,\,\mathsf{char}(\k)=0$ and $M$ be a finitely generated $\D(X)$-module.
 We also choose a projective embedding $i:T^* X\hookrightarrow \P^N_\k$.
Noetherianity of $\D(X)$ implies that $M$ is the cokernel of a morphism of free finitely generated $\D(X)$-modules.
Therefore, there exists a finitely generated ring $R\subset \k$ such that variety $X$, embedding $i$  and module
$M$ have models $X_R,i_R,{M}_R$
over $\spec  R$. We assume that $X_R$ is a smooth affine variety over $\spec R$ and $M_R$ is a finitely presented $\D(X_R)$-module.

For any  prime $p$ we obtain a finitely generated module $M_R\otimes \Z/p\,\Z$ over noncommutative ring 
$\D({X}_R)\otimes \Z/ p\,\Z$ which is a finitely generated module over its center. In particular, we can consider
 $M_R\otimes \Z/p\,\Z$ as a module over the center (it is again finitely generated). Hence, for any prime $p$ and 
for any point $v\in \spec R_p$ (here $R_p:={R}\otimes\Z/p\,\Z$), 
 with the residue field $\k_v$,  we obtain a reduced subscheme over $\k_v$
$$\mathsf{supp}_{p,v}( M_R)\subset T^*X'_{R_p}\times_{\spec  {R}_p}\spec \k_v\,.$$
Here $X'_{R_p}$ is the Frobenius pullback of $X_{R_p}:=X\times_{\spec R}\spec R_p$.
Then one has the following easy result (the proof is omitted here).
\begin{prop}
The dimension of $\mathsf{supp}_{p,v}({ M}_R)$ coincides with the dimension of $\mathsf{supp}\,(M)$ for large enough prime (and all
points $v\in \spec R_p$).
\end{prop}
Projective embedding $i_R:T^*_{X_R}\hookrightarrow \P^N$ induces an embedding $i'_{R_p}:T^*_{X'_{R_p}}\hookrightarrow \P^N$.
Therefore, one can speak about the degree of $\mathsf{supp}_{p,v}({ M}_R)$ via the embedding $i'_{R_p}$, as the sum of degrees
  of closures in $\P^N_{\k_v}$ of generic points of top-dimensional components of $\mathsf{supp}_{p,v}({ M}_R)$. It seems that following holds:
\begin{conj} {\rm The degree of $\mathsf{supp}_{p,v}( M_R)$ is bounded above by $const\cdot p^r$ where
$r=\dim \mathsf{supp}\,(M)- \dim X$.}
\end{conj}
We have a good evidence for this conjecture in the case of cyclic $\D(X)$-modules  of the form $M=\D(X)/\D(X)\cdot P$
where $P\ne 0$ is a non-zero differential operator, see section 3.2 for the special case $X=\A^1_\k$ and section 5.1 for $X=\A^n_\k,\,\,n>1$.

In the case of holonomic $\D$-modules the conjecture implies that the degree is uniformly bounded.
Also, one can expect an analog of Gabber theorem:
\begin{conj} {\rm In the above notation, for holonomic $M$ the support $\mathsf{supp}_{p,v}{ M}_R$ is Lagrangian for sufficiently large $p$ and any $v$.}
\end{conj}
Let us assume the above conjecture, and let $L_i$ be the generic point of an irreducible component 
 of $\mathsf{supp}_{p,v}({ M}_R)$. Then the length of $${ M}_v:={ M}_R\otimes_R \k_v$$ at $L_i$ is 
divisible by $p^{\dim X}$.
 Hence, we have an effective algebraic cycle on $T^*X_v'$ where $X_v':={ X}_{R_p}'\times_{\spec R_p}\spec \k_v $ given by
$$\mathsf{supp}_{p,v}^{\op{num}}({M}_R):=\sum_i \frac{\op{length}_{L_i} { M}_v}{p^{\dim X}} [L_i]\,.$$

\subsection{Arithmetic support}

Let us use the notation from  section 2.1, and assume that module $M$ is  holonomic. 
Let us assume\footnote{My student Thomas Bitoun informed me recently that he proved
conjectures 1 and 2. All the further conjectures made in the present paper seem to be completely out of reach.}  conjectures 1 and 2. Then  there exists prime $p_0$
 such that for any $p\ge p_0$ and any $v\in\spec R_p$ we have an effective Lagrangian cycle $\mathsf{supp}_{p,v}^{\op{num}}({M}_R)$ in $T^*X'_v$. Let us replace $R$ by its localization obtained by inverting all primes $<p_0$.

\begin{dfn}  The arithmetic support $\mathsf{supp}_{arith}^{\op{num}}(M)$ of  $M$ is an element of $ELC(T^*X)_\infty^{tw}$ uniquely specified by the condition
 that for a model ${ M}_R,{ X}_R$ over a finitely generated ring $R\subset \k$ as above, it is given by 
 the collection of Lagrangian cycles $\mathsf{supp}_{p,v}^{\op{num}}({ M}_R)$.
\end{dfn}

It is easy to see that the definition is consistent, i.e., that  for any model the collection  of cycles comes from a collection of constructible maps
 for all sufficiently large primes, and that  models form a filtered system. 

The main advantage of the arithmetic support is that it gives a more elaborate
 signature of a holonomic module, and the Lagrangian cycle is no longer conical in general.
 In the next subsection we will give explicit non-conical examples of arithmetic supports.
 
 The usual support  describes only the limiting behavior at infinity (along fibers of the cotangent bundle $T^*X\to X$) 
 of the arithmetic one. More precisely, there is a projection $ELC(T^*Y)\to ELC(T^*Y)$
  which associates with every Lagrangian cycle its limit under the contraction by $\lambda\in 
  \G_m,\,\,\lambda\to 0$. This limit is automatically conical.
  We expect that the conical limit of the arithmetic support  coincides with the pullback by the Frobenius of the usual support.
    The same can be said about the Bernstein filtration and the corresponding support
     in the case $X=\A^n_\k$.
     
     The arithmetic support is  covariant with respect to automorphisms of $X$ (and also contravariant for \'etale maps), as well as under symplectic affine transformations of $\A^{2n}_\k=T^*\A^n_\k$ in the case $X=\A^n_\k$. Hence, we see that it naturally generalizes
      two classical types of supports.

 Obviously, the arithmetic support behaves additively for extensions of $\D_X$-modules, hence it is sufficient to study it only for simple holonomic modules. Also, if $\mathsf{supp}_{p,v}^{\op{num}}({ M}_R)$ for infinitely many pairs $(p,v)$ with $p\to \infty$
 is just one irreducible Lagrangian subvariety taken with multiplicity one, then $M$ is simple.

Morally, we should think about the arithmetic support as about Lagrangian cycle defined over $\k_\infty$, this idea
 is elaborated further in section 3.1. The version with constructible maps presented here is a surrogate
 for the ``right'' version in section 3.1.

\subsection{Examples: $p$-curvature, exponents, fractional powers and Gauss-Manin connections}

Let $M$ be holonomic $\D_X$-module corresponding to a vector bundle $\cal E$
over $X/\k$ with flat connection $\nabla$. Let us choose a model of $(X,\cal E,\nabla)$ over a finitely generated ring 
$R\subset \k$. Then for each prime $p$ we obtain a bundle ${\cal E}_p$ with flat connection over a smooth scheme $X_p/ \spec R_p$
in characteristic $p$, where $R_p:=R\otimes \Z/p\,\Z$.
 The $p$-curvature of such a connection is a $p$-linear map
$$T_{X_p/\spec R_p}\to {{\cal E}nd}\,{\cal E}_p,\,\,\,\xi\mapsto \left(\nabla_\xi\right)^p-\nabla_{\xi^p} \,\,. $$
Moreover, the image of this map consists of commuting operators, hence we can interpret
 $p$-curvature as the Higgs bundle structure on ${\cal E}$. More precisely,
 it is a coherent sheaf ${\cal E}_{Higgs}$ on $T^*X'_p$ where $X'_p/\spec R_p$ is the pullback of $X_p$ under the Frobenius map
 ${\rm Fr}_{R_p,p}: \spec R_p\to \spec R_p$, together with an isomorphism of coherent sheaves on $X'_p$:
   $$(\op{pr}_{T^* X'_p\to X'_p})_* {\cal E}_{Higgs}\simeq (\op{pr}_{X'_p\to X})^* {\cal E}\,.$$
It follows directly from definitions that the arithmetic support 
 of $M$ (at prime $p$) is the same as the support of ${\cal E}_{Higgs}$ (compare with \cite{OgusVologodsky}).

There are two cases when one can easily calculate the arithmetic support.
First, for any $F\in \O(X)$ we have an associated holonomic $\D_X$-module given by the
 trivial line bundle $\O_X$ endowed with the flat connection
$$\nabla=d+(dF\wedge\cdot)\,.$$
One can think about this $\D_X$-module as $\exp(F)\cdot \O_X$.
We claim that the corresponding arithmetic support is the pullback by the universal Frobenius $\op{Fr}_{\k_\infty}$ of the graph
 of differential form $-dF$. This follows easily from the identity
 $$\left(\frac{d}{dx}+\frac{dG}{dx}\right)^p=\left(\frac{d}{dx}\right)^p+\left(\frac{dG}{dx}\right)^p$$
 which is held in $\D(R[x])$ for any ring $R$ over $\Z/p\,\Z$ and any element $G\in R[x]$ (see proposition 3 in \cite{KaKo}).
In particular, for the case $X=\A^n_\k$ and $F$ polynomial of degree $\le 2$, the support is the pullback by the Frobenius of
 the affine Lagrangian subspace in $\A^{2n}_\k=T^*X$ corresponding to $F$.
 
One can also calculate the arithmetic support for the connection on the trivial bundle
 corresponding to a closed but not exact 1-form. For example,
 for $X=\spec \k[x,x^{-1}]$ and for the connection on ${\cal E}:=\O_X$ given by $1$-form $\lambda dx/x$ for some $\lambda\in \k$,
 the arithmetic support is the curve in $T^*X$ given by the equation
$$x_{(p)}y_{(p)}=\lambda ^p-\lambda {\pmod p}\,.$$
Here $x_{(p)},y_{(p)}$ are coordinates on $T^*X'_p\subset T^*\A^1_\k$ with symplectic form $dx_{(p)}\wedge dy_{(p)}$.
In the case $\lambda\in \Q$ the expression $(\lambda^p-\lambda) \pmod p$ vanishes for all sufficiently large $p$, hence the arithmetic support of this $\D_X$-module is just the zero section of $T^*X'_p$. If $\lambda\notin \Q$ then the 
arithmetic support is not equal to the zero section, as follows from Chebotarev density theorem in the case
 when $\lambda$ is algebraic, and by elementary reasons when $\lambda$ is transcendental.

Finally, let $M$ be a holonomic $\D_X$-module corresponding to the vector bundle $\cal E$ on $X$ endowed with a flat connection $\nabla$ of Gauss-Manin type (for variations of pure motives).
 This means that $({\cal E},\nabla)$ is a subquotient of the natural connection of the bundle of de
 Rham cohomology of fibers of a smooth projective morphism $Y\to X$. Then by a classical result of N.~Katz (see \cite{Katz}) the
 $p$-curvature is nilpotent. Hence the support is the zero section of the cotangent bundle, taken with the multiplicity
 equal to $\op{rank} {\cal E}$.

\section{One-dimensional case}

\subsection{From higher-dimensional case to $\A_\k^1$}

Here we sketch a geometric construction which reduces the study of arithmetic supports in higher dimensions to the case of $\A^1_\k$. 

First of all, the arithmetic support is compatible with localization, hence we can assume that we consider 
holonomic $\D(X)$-modules for a smooth {\it affine } variety $X$. Let us choose a closed embedding $j:X\hookrightarrow
 \A^N_\k$ for some $N\in \Z_{\ge 0}$. Then any holonomic $\D$-module $M$ gives a holonomic $\D$-module $j_* M$
 on $\A^N_k$. The behavior of 
arithmetic supports for embeddings is the natural one, given by the Lagrangian correspondence in $T^*X\times T^*\A^N_\k$ equal to the conormal bundle to $\op{graph}(j)$.

 Now, we describe a way to put a structure of a reduced ind-scheme on the ind-constructible set $ELC(\A^{2N}_\k)$, similar to 
 the Chow scheme (see \cite{Angeniol}) parametrizing algebraic cycles of given dimension and degree in a smooth
projective scheme. Here $\A^{2N}_\k=T^*\A^N_\k$ is considered as a symplectic manifold. 

Let us consider the variety $V/\k$ parametrizing triples $(B,b_1,b_2)$ where $B\subset \A^{2N}_\k$ is a coisotropic
 affine subspace in $\A^{2N}_\k$, of dimension $N+1$, and $b_1,b_2$ are two points in the symplectic affine plane $\widetilde{B}\simeq \A^2$ which 
is obtained by factorization of $B$ along the kernel of the natural Poisson structure on $B$.
We claim that every effective Lagrangian cycle $C=\sum_{i\in I} m_i L_i$ gives a non-zero rational function $\phi_C$ on $V$.

Indeed, for generic $B$ the intersection of all subvarieties $L_i$ with $B$ is one-dimensional and transversal
at the generic point of every component. Moreover, its projection to $\widetilde{B}$ is a plane curve. Hence, taking the sum over $i$ 
we obtain an effective divisor $C_B$ on plane $\widetilde{B}$ (a collection of curves with positive multiplicities).
 There exists a unique up to scalar non-zero polynomial $F_B$ on $\widetilde{B}$ whose divisor of zeroes is $C_B$.
Then, for generic $b_1,b_2\in \widetilde{B}$ the ratio $F_B(b_1)/F_B(b_2)$ is canonically defined and is not zero. We set
$$\phi_C(B,b_1,b_2):=F_B(b_1)/F_B(b_2)$$ for generic $(B,b_1,b_2)$.

Let $R/\k$ be a finitely generated algebra {\it without nilpotents}.
We define a family over $\spec{R}$  of effective Lagrangian cycles on $\A^{2N}_\k$ to be 
 a pair $(U,\phi)$ where $U\subset \spec R\times_{\spec k} V$ is a Zariski open subset which is dominant
 over $\spec R$, and $\phi\in \O(U),\phi\ne 0$ is a function on $U$ such that for every point $x\in \spec R$
 the restriction of $\phi$ to the fiber over $x$ coincides with the restriction of a rational function associated
 with an effective Lagrangian cycle on $\A^{2N}_{\k_x}$. We identify two pairs $(U,\phi)$ and $(U',\phi')$ if and only if
 $$\phi_{|U\cap U'}=\phi'_{|U\cap U'}\,.$$

Thus, we have defined $ELC(\A^{2N}_\k)$ as a set-valued functor on finitely generated rings without nilpotents.
One can check that this gives a structure of a reduced ind-scheme.
The above definition (at least of a functor) also work without the assumption that the ground field $\k$ has zero characteristic, and it 
generalizes immediately to symplectic manifolds over arbitrary base.

Similar constructions can be performed for holonomic $\D(\A^n_\k)$-modules. The analog of intersection with $B$ and projection to $B'$ is given by the functor from holonomic $\D(\A^n_\k)$-modules to holonomic $\D(\A^1_\k)$-modules determined by the kernel
 corresponding to an affine Lagrangian subspace in $T^*(\A^n_\k\times \A^1_\k)$.
 One expects that the arithmetic supports for such functors behave as is prescribed by the geometric construction from above.  Hence, one get a reduction of the problem of the description of the arithmetic support to the case $X=\A^1_\k$.

Moreover, it seems that there should exist an enhanced definition of the arithmetic support of a holonomic $\D_X$-module $M$
 as  a $\k_\infty$-point of the ind-scheme ${\op{Fr}}^*_{\k_\infty}ELC(T^*X)$.

\subsection{Arithmetic support of a cyclic module}

Let us consider the case $X=\A^1_\k$, $\mathsf{char}(\k)=0$. The algebra $\D(X)$ is the first Weyl algebra $A_{1,\k}$, we denote its generators by $\xh=x$ and ${\yh}=d/dx$.
 We consider   holonomic $\D(\A^1_\k)$-module which is a non-trivial cyclic module
   $\D(\A^1_\k)/\D(\A^1_\k)\cdot P$, where
$$P=\sum_{i+j\le N} a_{ij} x^i (d/dx)^j$$
is a non-zero differential operator on $X$.
Let us fix a  finitely generated subring $R\subset \k$ containing all coefficients
   $a_{ij}$. The center of $A_{1,R_p}$ (recall $R_p:=R\otimes \Z/p\,\Z$) is the polynomial algebra
  $R_p[\xh^p,{\hat {y}}^p]$. We extend it by adding central variables $\tilde{x},\tilde{y}$ satisfying
  $$\tilde{x}^p=\xh^p,\,\,\,\tilde{y}^p={\hat {y}}^p\,\,.$$
 The resulting extension of $A_{1,R_p}$ is isomorphic to the matrix algebra 
$${\rm Mat}(p\times p,
   R_p[\tilde{x},\tilde{y}])\,.$$ Indeed, this extension is the algebra over $R_p[\tilde{x},\tilde{y}]$ generated
 by two elements $\xh,\yh$ satisfying the relations 
   $$[\yh,\xh]=1,\,\,\xh^p=\tilde{x}^p,\,\,\yh^p=\tilde{y}^p\,.$$
Shifted generators $(\xh-\tilde{x},\yh-\tilde{y})$ satisfy the same relations as the operators $x$ and $d/dx$
 in the truncated polynomial ring $\Z/p\,\Z\,[x]/(x^p)$:
$$x^p=0, (d/dx)^p=0, [d/dx,x]=1\,. $$ 
For example, for $p=5$ the corresponding matrices are
$$X_p=\left(\begin{array}{ccccc} 0 & 0 & 0 & 0 & 0\\
                           1 & 0 & 0 & 0 & 0\\
                           0 & 1 & 0 & 0 & 0\\
                           0 & 0 & 1 & 0 & 0\\
                           0 & 0 & 0 & 1 & 0
      \end{array}\right),  \,\,
 Y_p=\left(\begin{array}{ccccc} 0 & 1 & 0 & 0 & 0\\
                           0 & 0 & 2 & 0 & 0\\
                           0 & 0 & 0 & 3 & 0\\
                           0 & 0 & 0 & 0 & 4\\
                           0 & 0 & 0 & 0 & 0
      \end{array}\right)\,\,.$$

Then one get the following description of the arithmetic support of $M=\D(\A^1_\k)/\D(\A^1_\k)\cdot P$.
Namely, for a given prime $p$ let us consider the polynomial 
$$ \widetilde{D}_{P}^{(p)}:= \det\left(\sum_{i+j\le N} a_{ij} \left( X_p+\tilde{x}\cdot \bf{1}_p\right)^i\cdot 
\left( Y_p+\tilde{y}\cdot \bf{1}_p\right)^j\right)\in R_p[\tilde{x},\tilde{y}]$$
where $\bf{1}_p$ is the identity matrix of size $p\times p$.
Obviously $\widetilde{D}_P^{(p)}$  has degree $\le N\cdot p$ in $\tilde{x},\tilde{y}$.
 We claim that it is in fact a polynomial of degree $\le N$ in $\tilde{x}^p,\tilde{y}^p$.
 This can be seen by general reasons, and as well by a direct check. Namely, the property of a polynomial
 in characteristic $p$ to depend only on $p$-th powers of variables is equivalent to the vanishing of its partial derivatives:
$$\frac{\partial}{\partial \tilde{x}} \widetilde{D}_P^{(p)}=\frac{\partial}{\partial \tilde{y}} \widetilde{D}_P^{(p)}=0\,.$$
The vanishing of, say, the derivative with respect to $\tilde{x}$ can be proved as follows. Taking this derivative is equivalent to taking the derivative of the determinant of the matrix from above under the infinitesimal
 conjugation by ${\bf{1}_p}+\epsilon Y_p$ where $\epsilon$ is a small parameter, $\epsilon^2=0$.
 The invariance of the determinants under the conjugation proves the result.

Therefore, we can write
  $$ \widetilde{D}_{P}^{(p)}(\tilde{x},\tilde{y})=D_P^{(p)}(\tilde{x}^p,\tilde{y}^p)$$
where $D_P^{(p)}$ is a polynomial in two variables of degree $\le N$:
$$D_P^{(p)}\in R_p[\tilde{x}^p,\tilde{y}^p]=R_p[\xh^p,\yh^p]=\op{Center}(A_{1,R_p})\,.$$

\begin{prop} The arithmetic support of module $M$ at pair $(p,v)$ where $p$ is a prime and $v\in\spec R_p$,
 is the effective one cycle on plane $\A^2_{\k_v}=\op{Fr}^*_{\k_v} \A^2_{\k_v}$ given as the divisor of zeroes of the image of polynomial $D_P^{(p)}$
 in $\k_v[\tilde{x}^p,\tilde{y}^p]$.
\end{prop}

This proposition follows directly from the definitions, and from the obvious equivalence
$$\op{Mat}(p\times p, \k)/\op{Mat}(p\times p, \k)\cdot T\ne 0\,\,\,\Longleftrightarrow
\,\,\,\op{det} (T)=0$$
for any matrix $T\in \op{Mat}(p\times p, \k)$ and any field $\k$.

The free term of the polynomial $D_P^{(p)}$ is
$$\op{det}_p(P):=\det\left(\sum_{i+j\le N} a_{ij} X_p^i\cdot 
Y_p^j\right)\in R_p\,\,.$$
This expression we will call the {\it $p$-determinant} of a polynomial differential operator in one variable.
One can treat coefficients $(a_{ij})_{i+j\le N}$ as independent variables, hence we have a {\it universal} 
$p$-determinant
$$\op{det}_p^{\le N}\in \Z/p\,\Z[(a_{ij})_{i+j\le N}]$$
which is a homogeneous polynomial of degree $p$ in $\frac{(N+1)(N+2)}{2}$ variables with coefficients in $\Z/p\,\Z$.

The calculation of other coefficients of $D_P^{(p)}$ can be reduced  to the calculation of finitely many 
$p$-determinants of  differential operators. Indeed,  any  polynomial
 of a bounded degree can be reconstructed by the Lagrange interpolation formula from  its values at finitely many points.

\subsection{Determinant formulas}

One can calculate $p$-determinants effectively using a well-known formula. The algorithm runs very fast, linearly in prime $p$.
  Notice that the matrix $$M_p:=\sum_{i+j\le N}a_{ij}X_p^i\cdot 
Y_p^j$$ contains non-zero terms only at distance at most $N$ from the main diagonal. We are interested in its determinant for $p\gg N$.

Consider the general situation: let $M$ be a square matrix of size $L\times L$ (with coefficients in a commutative ring)
such that $M_{ij}=0$ if $|i-j|>N$ for some $N<L/2$. For every integer $i\in [N+1,L]$ denote $A^{(i)}$ the square matrix of size $2N\times 2N$
 given by 
$$A^{(i)}_{j_1,j_2}=\left\{\begin{array}{cl}  M_{i-N,i}\hfill& \ \mbox {if } \  j_1=j_2-1 \\
   -M_{i-N,i-2N+j_2-1} &\mbox { if } \  j_1=2N \hfill \\
  0 \hfill &\mbox{ otherwise }  \hfill \end{array}\right.$$
Here we set $M_{ij}:=0$ for $j\le 0$.
Also introduce rectangular matrices 
$$ B\in \op{Mat}(N\times 2N),\,\,B_{j_1,j_2}=M_{j_1+L-N, j_2+L-2N}\,\,,1\le j_1\le N,\,\,1\le j_2\le 2N\,,  $$
$$B'\in \op{Mat}(2N\times N),\,\, B'_{j_1,j_2}=\delta_{j_1,j_2-N}\,\,,1\le j_1\le 2N,\,\,1\le j_2\le N\,.$$
\begin{prop} In the above notation one has
$$ \left(\prod_{i={N+1}}^L M_{i-N,i} \right)^{N-1}\cdot\det(M) = \pm\det\left( B\cdot\left( A^{(L)} A^{(L-1)}\cdots A^{(N+1)}\right) \cdot B'\right)\,.$$ 
\end{prop}
{\it The idea of the proof.} Suppose that the matrix $M$ is degenerate and all elements $M_{i-N,i}$ are non-zero for $i\in [N+1,L]$.
Hence the left hand side of the above identity vanishes, and we want to prove that the right hand side vanishes too. 
 Let us  consider the sequence
$$(v_j)_{j=1,N+L}:=(\underbrace{0,\dots,0}_{N\,\,\mathrm{times}},x_1,\dots, x_L)      $$
where $(x_1,\dots,x_L)$ is a non-zero vector in the kernel of $M$.
Then its subsequences of length $2N$
$$U^{(j)}:=(v_j,v_{j+1},\dots ,v_{j+2L-1}),\,\,\,1\le j\le L-N+1$$ 
satisfy the relations
$$\begin{array}{l} U^{(1)}\in \op{Im}(B') \\
 U^{(2)}=M_{1,N+1}^{-1} \cdot A^{(N+1)} \,U^{(1)}\\
 U^{(3)}=M_{2,N+2}^{-1}\cdot A^{(N+2)} U^{(2)}\\
\dots\\
U^{(L-N+1)}=M_{L-N,L}^{-1}\cdot A^{(L)} \,U^{(L-N)}\\
0=B \,U^{(L-N)}\,.
\end{array}$$
Hence we conclude that 
$$\det \left( B\cdot\left( A^{(L)} A^{(L-1)}\cdots A^{(N+1)}\right) \cdot B'\right)=0\,.$$
\hfill $\Box$

The above proposition allows to calculate $p$-determinants up to a simple factor which has the  form 
$$\mbox{ prime } p\mapsto \prod_{j=1}^{p-N} f(j)\pmod{p}\,,$$
where $f=f(x)\in R[x]$ is as polynomial in one variable with coefficients in a finitely generated ring $R\subset \k$.
Such factors vanish sometimes (e.g. when $f$ has a root in $\Q$), in this case one modify it by replacing $f$ by $f+c$ where $c$ is a new independent 
constant. 

We see that $p$-determinants (up to factors discussed above) belong to the following class of expressions 
$$\mbox{ prime }p\mapsto \op{Tr}(F(1)\cdot F(2)\cdot\ldots\cdot F(p-k)\cdot G)\pmod{p} ,\,\,\mbox{ if } p\ge k$$
where $F\in \op{Mat}(K\times K,R)[x], \,G\in \op{Mat}(K\times K,R)$ for some $K,k\in \Z_{\ge 1}$.

\subsection{Logarithmic families of planar curves}

Let $\k$ be an algebraically closed field of characteristic zero. The set of planar curves $Curves_{\A^2_\k}$ understood as effective divisors in $\A^2_\k$, is the same as the quotient of the set of non-zero polynomials $P\in \k[x,y]\setminus\{0\}$ modulo multiplicative constant. 
Hence it carries a natural structure of an ind-scheme, it is an infinite-dimensional
 projective space
 $$Curves_{\A^2_\k}=\P^\infty(\k)=\P(\k[x,y])=\lim_{\longrightarrow}\P^{\frac{(d+1)(d+2)}{2}-1}(\k)\,.$$
 The group  $\op{Aut}\A^2_\k$ acts by automorphisms of this ind-scheme.

Our goal here is to introduce an equivalence relation on $Curves_{\A^2_\k}$ invariant under $\op{Aut}\A^2_\k$ 
 such that all equivalence classes will be sets of $\k$-points of constructible sets, and for any 
 cyclic $\D(\A^1_\k)$-module $M$ its arithmetic support will belong to one such an equivalence class.
   These equivalence classes we will call {\it logarithmic families} because they have a characterization
 in terms logarithmic divergence of certain integrals.

First, we introduce certain set $J^\infty_\k$ associated with $\A^2_\k$.
 It can be thought as truncated jets of algebraic curves in $\P^2_\k$ intersecting
  the projective line at infinity $\P^1_k=\P^2_\k\setminus \A^2_\k$.
 The truncation  means that we do not specify the terms in Puiseux series
 which change the germ of the curve by another germ such that the area with respect to 2-form
  $dx\wedge dy$ in a segment bounded by germs is finite.
  Here is the precise definition for the special case of germs intersecting $\P^1_\k$
  at the infinite point on the $x$-axis given by $y=0$ in $\A^2_\k=\spec \k[x,y]$.
   The germ is given by an integer $d\ge 1$ and a sequence of numbers $a_i\in \k, \,-d<i<d$ such
   that $$g.c.d.\left(\{d\}\cup\{i\,|\,a_i\ne 0\}\right)=1\,.$$
   We identify such data factorizing by the free action of the group $\mu_d$ of roots 1 of order $d$
   $$\xi\in \k,\,\xi^d=1\mbox { acts as } a_i\mapsto \xi^i a_i\,.$$
   The interpretation of the pair $d,(a_i)_{-d<i<d}$ is as the truncated germ of a curve
   $$x=x(t)=t^d,\,\,\,y=y(t)=\sum_{1-d}^{d-1} a_i t^i +O(t^{-d}),\,\,t\to \infty \, .$$
  Alternatively, one can write $y$ as a Puiseux series in $x$:
     $$y=a_{d-1} x^{\frac{d-1}{d}}+\dots +a_{1-d}x^{-\frac{d-1}{d}}+O(1/x)\,\,.$$
     Conversely, any Puiseux series 
     $$y\sim\sum_{\lambda\in \Q,\,\lambda<1} c_\lambda x^\lambda $$
     which is $\sim o(x)$ as $x\to \infty$, has a unique representation as above \break $\pmod{O(1/x)}$.
      Namely, we define $d$ as the minimal integer $\ge 1$ such that
       $\lambda\in (1/d)\cdot \Z$ for all $\lambda\in (-1,1)$ such that $c_\lambda\ne 0$,
        and set $a_i:=c_{i/d}$.
        
        Acting by the group $GL(2,\k)$ we obtain the description of the whole set $J^\infty_\k$. There is a $\Z_{\ge 1}$-valued function $deg$ on $J^\infty_\k$, with the value equal to $d$ in the above notation.
        
        Any planar curve $C\in Curves_{\A^2_\k}$ gives a function $\nu_C:J^\infty_\k\to \Z_{\ge 0}$ with finite support. Namely, we count with multiplicities all the germs of $C$ intersecting $\P^1_\k$ at infinity. The degree of $C$ coincides with the sum over $J^\infty_\k$ of the product of $\nu_C$ with $deg$.
        
        We define a logarithmic family to be the set of all curves $C$  with a given function  $\nu_C$. 

Here are examples in small degrees.
        First, we have the logarithmic family consisting of the empty curve $C$ with $\nu_C=0$. Next, we have one-point logarithmic families each of which consists
         of a line in $\A^2_\k$.  The simplest non-trivial example is the 
         family of hyperbolas $xy=t$ where $t\in \k$, including the degenerate case $t=0$.
         
        It is easy to see that  there is a natural action of $\op{Aut}\A^2_\k$ 
        on $J^\infty_\k$, and the decomposition by logarithmic families is
       $\op{Aut}\A^2_\k$-equivariant. An intrinsic definition of $J^\infty_\k$ is 
       as  the inductive limit of the set of divisors where the volume form $dx\wedge dy$ has pole of order one,
        over the partially ordered set of smooth compactifications of $\A^2$ on which the
        volume form does not vanish at infinity (i.e., it is a Poisson compactification, compare with \cite{Kontsevich1}). For any curve $C\in Curves_{\A^2_\k}$ there exists a Poisson compactification of $\A^2_\k$ such that $C$ intersects only those divisors at infinity
        where the form has logarithmic pole. Also, If $C_1$ and $C_2$ are two curves
        such that $\nu_{C_1}\cdot \nu_{C_2}=0$ (i.e., functions $\nu_{C_1}$ and $\nu_{C_2}$
         have disjoint support), then the intersection $C_1\cap C_2$ is finite and the intersection number $[C_1]\cap[C_2]$ can be determined entirely in terms of functions
   $\nu_{C_1}$ and $\nu_{C_2}$.  
   
   For a non-algebraically closed field $\k$ of zero characteristic we define logarithmic
       families using the embedding $\k\hookrightarrow \overline{\k}$ to the algebraic closure.
  Also, logarithmic families for curves of a given bounded degree can be defined for positive characteristic if it is large enough.

       There is an alternative meaning of $J^\infty_\k$ in terms of singularities of holonomic 
       $\D(\A_{\k})$-modules. Let $M$ be such a module. We will associate with $M$ a
        $\Z_{\ge 0}$-valued function $\nu_M$ on $J^\infty_\k$ with finite support.
        
        First,  there exists a finite 
       set $S\subset \k$ such that $M_{|\A^1_\k\setminus S}$ is a vector bundle $\cal E$ with 
       connection $\nabla$. It can have irregular singularities at $S$ and at $\infty$.
       
       It is well-known (see e.g. \cite{Malgrange}) that the category of bundles with connections over the field of Laurent series $\k((z))$ is decomposed into the direct sum of blocks corresponding to  Puiseux polynomials in negative powers of z:
       $$F(z)=\sum_{\lambda\in \Q_{<0}} b_\lambda z^\lambda,\,\,b_\lambda\in \k,\,\,b_\lambda=0 \mbox { for almost all }\lambda\,,$$
       defined modulo the action of $\mu_d$ where $d$ is $l.c.m.$ of all denominators
        of $\lambda$ with $b_\lambda\ne 0$. The basic $\D$-module in such a block is
         $\exp(F)\cdot\k((z))$.
         These blocks for $F\ne 0$ correspond exactly to elements in $J^\infty_\k$
         corresponding to truncated germs intersecting the divisor at infinity at point
          $(0,\infty)$ in the compactification 
          $$\A^2_\k=\A^1_\k\times\A^1_\k\subset\P^1_\k\times\P^1_\k\,,$$
           with the exception of the germ at infinity of the vertical line $z=0$.
           The correspondence is given by
            $$F\mapsto \mbox { germ of the curve } (z,F'(z))\,\mbox{ at } z \to 0\,.$$
            
            In this manner we will associate to any holonomic $\D(\A^1_{\k})$-module $M$ multiplicities $\nu_M$
             at all points of $J^\infty_\k$ except truncated germs of lines
             $x=x_0$ for $x_0\in S$. In order to get multiplicities at these points
              one can apply an automorphism of the Weyl algebra $A_{1,\k}$ corresponding to
               a non-trivial matrix in $SL(2,\k)$, e.g. the Fourier transform.  
    We notice also that the complicated formulas from \cite{Malgrange} relating irregular singularities of $M$ with the ones of its Fourier transform, translate just to the  action of matrix $\left(\begin{array}{cc} 0 & 1\\-1 &0 \end{array}\right)$
 on $J^\infty_\k$.

           A holonomic $\D(\A^1_{\k})$-module $M$ has only regular singularities if and only if the multiplicities 
$\nu_M$
           vanish  at all points of $J^\infty_\k$ except the germs of lines given by equations
$$ y=0,\,\mbox{ or }\,x=x_0\mbox{ for some }\,x_0\in \k\,.$$

       \begin{thm} For any non-zero differential operator $P\in A_{1,\k}=\D(\A^1_\k)$
        with coefficients in a finitely generated ring $R\subset \k,\,\,\mathsf{char}(\k)=0$,
         for all sufficiently large $p$ and for any point $v\in \spec R$ over $p$,
          the multiplicities $\nu_M$ at point $v$ where $M=A_{1,\k}/A_{1,\k}\cdot P$ coincides
           with the (pullback by Frobenius of) multiplicities of the curve 
              given by the equation $D^{(p)}_P=0$ at point $v$.
       \end{thm}

{\it The idea of the proof}. First of all, it is easy to identify
contributions of germs of the line $y=0$. Namely, in the case of a curve given by equation
 $H(x,y)=0,\,H=\sum_{i,j} H_{ij} x^i y^j\ne 0\in \k[x,y]$, this multiplicity can be read from the Newton polygon of $P$.
   Namely, the multiplicity is equal to
   $$\max\{j\,|\, H_{ij}\ne 0,\,\forall {(i',j')} \,H_{i'j'}\ne 0\Longrightarrow
     (i-j)\ge (i'-j')\} \, .$$
  A similar description works for 
 cyclic $\D_{\A^1_\k}$-modules.
 For the multiplicities at other points of $J^\infty_\k$ one can apply automorphisms
  of the Weyl algebra, and also take the tensor product with $\D$-modules corresponding to
   exponents of Puiseux polynomials. \hfill $\Box$

Finally, one can show  that for two holonomic $\D_{\A^1_\k}$-modules $M_1,M_2$ such that
supports of $\nu_{M_1}$ and $\nu_{M_2}$ are disjoint, there is no non-trivial homomorphisms
 from $M_1$ to $M_2$, and the dimension of $\op{Ext}^1(M_1,M_2)$ coincides
  with the intersection number of the corresponding curves.

\section{Correspondence between classical and quantum families}

\subsection{Rough picture for Lagrangian cycles}

We expect that in the case $\dim X>1$ also there exists a notion of a logarithmic family of effective Lagrangian cycles
in $T^*X$, and the arithmetic support should always belong to such a family.
 In the special case when a Lagrangian cycle is a {\it smooth} closed Lagrangian variety $L\subset T^*X$
 (taken with multiplicity one) we expect  a more clearer picture of what is the logarithmic family:
 \begin{dfn} A smooth logarithmic family of smooth Lagrangian subvarieties in $T^*X$ 
   is a pair $(S,{\cal L})$ where $S$ is a smooth variety over $\k$ and
    ${\cal L}\subset T^*X \times S$ is a smooth closed submanifold such that its projection
     to $S$ is smooth, all fibers ${\cal L}_s,\,s\in S$ are Lagrangian, and the following property holds.
     For any $s\in S$ the natural map
     $$T_s S\to \Gamma({\cal L}_s,(T_X)_{|{\cal L}_s}/T_{{\cal L}_s})=\Gamma({\cal L}_s,
     T^*_{{\cal L}_s})$$
     identifies $T_s S$ with the space of 1-forms on ${\cal L}_s$ with logarithmic singularities\footnote{All such forms are automatically closed.}.
 \end{dfn}
 \begin{conj} {\rm For a smooth closed Lagrangian $L\subset T^*X$ there exists a smooth logarithmic family $(S,{\cal L})$ with base point $s_0\in S$ such that ${\cal L}_{s_0}=L$. Also, any two such
  families coincide with each other in the vicinity of $s_0$.}
 \end{conj}
 
\subsection{Isosingular families of holonomic $\D$-modules}

By analogy with the geometry of logarithmic families, we expect that something similar
should happen for holonomic $\D_X$-modules as well . We say (to a first approximation) that two holonomic modules
 are isosingular if and only if the corresponding arithmetic supports belong to the same logarithmic family.
In the case $X=\A^1$ the precise definition is the coincidence of multiplicities, as explained in section 3.4.

Namely, for any holonomic module $M$ we expect that there exists a natural moduli stack $Mod_M$ parametrizing holonomic $\D_X$-modules
 which looks locally like the quotient of a scheme of finite type by  a group whose connected component of identity is a finite-dimensional affine algebraic group.
  The tangent complex of $Mod_M$ at the base point corresponding to $M$ should have cohomology
    $\op{Hom}(M,M)$ in degree $-1$ and $\op{Ext}^1(M,M)$ in degree $0$. Connected components
     of stacks $Mod_M$ we will call {\it isosingular} families of holonomic $\D_X$-modules.
    
    One of motivations is that the abelian category of holonomic $\D_X$-modules has the following finiteness property: for any two objects $M_1, M_2$ we have
    $$\dim \op{Hom}(M_1,M_2)<\infty,\,\,\,\dim\op{Ext}^1(M_1,M_2)<\infty\,.$$
    The same property is shared by the category of finite-dimensional modules over a finitely generated associative algebra $A/\k$.  In the latter case we have naturally defined moduli stacks
     of objects.
     
     Also, in the case $X=\A^1_\k$  any holonomic $\D_X$-module $M$ belongs to the abelian category
     consisting of all holonomic modules $M'$ such that 
     $$\mathsf{supp} \,\nu_{M'}\subset \mathsf{supp} \,\nu_M\,.$$
     In the case $\k=\C$ using the Riemann-Hilbert correspondence (for irregular singularities)
     we can identify the above abelian category with the category of finite-dimensional
     representations of a finitely generated algebra. Hence, we get moduli stacks of the
     form described above, but  with a ``wrong'' algebraic (but ``correct'' complex analytic)
     structure
      on the moduli stack.

     If the arithmetic support of a simple holonomic $\D_X$-module $M$ is a family of effective Lagrangian cycles such that the 
generic representative $L$ of such a family is a smooth connected non-empty closed subvariety with multiplicity one, then we expect that 
     $$\dim\op{Hom}(M,M)=\dim H^0(L)=1,\,\,\dim \op{Ext}^1(M,M)=\dim H^1_{dR}(L)\,.$$
     We see that the dimension of an isosingular family is greater than the dimension
      of the corresponding logarithmic family. The difference between two dimensions
      is the dimension of the Picard variety of any smooth compactification $\overline{L}$ of $L$.
      
      Informally speaking, holonomic $\D_X$-modules correspond to Lagrangian submanifolds in $T^*X$
      together with a line bundle on $\overline{L}$ (or something like that\footnote{In the picture with reduction modulo prime we get not only a Lagrangian submanifold but also a module over an Azumaya algebra, of minimal rank. Two such modules differ by a line bundle.}). More precisely, there should be the decomposition of the space of holonomic modules into isosingular families,
 and the decomposition of Lagrangian submanifolds with line bundles into logarithmic families, such that there is a canonical one-to-one correspondence between families (equivalence classes) of both kinds,
   and dimensions of the corresponding families coincide.

\subsection{Constants for the arithmetic support}

Let $\k$ be a field of characteristic zero, and consider the subring $P_\k\subset \k_\infty$ generated
 by $p$-determinants of differential operators in one variable, with coefficients in $\k$.
 Considerations from sections 3.1 and 3.2 lead to the following question. 
$$\mbox{\it What is the structure of }P_\k\,?$$
 The importance of the ring $P_\k$ is that the arithmetic supports should be  effective Lagrangian cycles parametrized by $\spec P_\k$.
 For countable $\k$ (e.g. for $\k=\Q$ or $\k=\overline{\Q}$) the ring $P_\k$ is a countable subring in uncountable ring $\k_\infty$, 
similar to the subalgebra of periods in $\C$ for algebraic varieties over $\Q$.

 The result from  section 3.3 indicates that $P_\k$ should be related with the differential Galois group of difference equations with coefficients in $\k$.

 Also, there are strong indications  that $P_\k$  should be related somehow to another group, the motivic Galois group of  $\k$. Indeed, for holonomic $\D_X$-module 
 corresponding to the flat connection on the trivial bundle given by a closed but not exact $1$-form, the support is controlled by the Cartier operator, which is related to the comparison of de Rham and crystalline cohomology of $X$ (in degree 1). 
Roughly speaking, one can expect that the motivic Galois group of $\k$ in de Rham realization acts on the category of holonomic $\D_X$-modules. 
The correspondence between holonomic modules and their arithmetic support could be related to torsors
 comparing de Rham and crystalline cohomology (and maybe Hodge realization which is the associated graded with respect to the Hodge filtration
 on de Rham cohomology). Fixed points of the motivic Galois group should correspond to motivic (or, more generally, motivic-exponential) holonomic  $\D_X$-modules. This fits well with the fact that supports of motivic and exponential $\D_X$-modules are defined over $\op{Fr}^*_{\k_\infty}(\k)\subset \k_\infty$.

For a given field $\k,\,\mathsf{char}(\k)=0$ we define {\it extended motivic-exponential} $\D_X$-modules on smooth algebraic varieties over $\k$
 as the minimal class which is closed under extensions, subquotients, pushforwards and pullbacks, and contain all $\D_X$-modules of type $\exp(F)\cdot \O_X$
 for $F\in \O(X)$. 
\begin{conj} {\rm The arithmetic support of a holonomic $\D_X$-module $M$ is the pullback by the universal Frobenius $\op{Fr}_{\k_\infty}$ of an effective Lagrangian cycle in $T^*X$ defined over $\k\subset \k_\infty$ if and only if $M$ is extended motivic-exponential.}
\end{conj}

\subsection{Isolated points}

 One can ask what are the ``most canonical'' holonomic modules and corresponding Lagrangian
 varieties\footnote{In a recent preprint \cite{Arinkin} a related but different question was studied for $X=\A^1_\k$.}. 
 \begin{conj} {\rm For any smooth closed connected Lagrangian subvariety $L$ in $T^*X$ over $\k=\C$ such that $H_1(L(\C),\Z)=0$ there exists a unique holonomic 
 $\D_X$-module $M=M_L$ with the arithmetic support equal to $L$ taken with multiplicity 1. Moreover,
  $\op{Ext^1}(M,M)=0$.}
 \end{conj}
 The reason  for the condition $H_1(L(\C),\Z)=0$ is that it guarantees that
  there is no non-trivial local system of rank 1 over $L$.
 This condition can be reformulated in  form which makes sense for arbitrary field $\k$ with $\mathsf{char}(\k)=0$:
$$H^1_{\acute{e}t}(L\times_{\spec \k}\spec{\overline{k}},\Z/l\Z)=0$$
for any prime $l$. 

 Presumably, one can weaken the condition on smoothness of $L$, e.g. it is definitely sufficient to assume that the codimension of singularities is $\ge 3$.

   \begin{conj} {\rm Any holonomic $\D_X$-module $M$ with $\op{Ext}^1(M,M)=0$ is of the extended  motivic-exponential type.}	
   \end{conj}

A corollary of conjecture 5 is one of conjectures discussed in \cite{KaKo}, which says
 that for any polynomial symplectomorphism $\phi$ of $\A^{2n}_\k$ there exists a canonically associated to it a bimodule $M_\phi$ over the Weyl algebra 
$A_{n,\k}$ giving Morita self-equivalence of the category of 
 $A_{n,\k}$-modules. Indeed, the graph of $\phi$ is a smooth Lagrangian subvariety in $\A^{4n}_\k$.
  Moreover, this subvariety is simply connected being isomorphic to $\A^{2n}_\k$.
  Hence it should give a canonical holonomic module $M_\phi$ over $A_{2n,\k}$ which can be interpreted as
   a bimodule over $A_{n,\k}$.

Also, conjecture 6 implies that $M_\phi$ is of the extended  motivic-exponential type.

\section{Relation with integrable systems}
\subsection{Arithmetic support of a non-holonomic module}

Let us consider a typical  non-holonomic module $M$ which is just a cyclic module $\D(X)/\D(X)\cdot P$  for a non-zero differential operator  $P\in \D(X)$. Here $X/\k$ is a smooth affine variety, of dimension $n>1$.
Then the support at prime $p$ is a hypersurface of degree which is bounded by $const\cdot p^{n-1}$.
This can be seen most easily in the case $X=\A^n_\k$. The consideration similar to one from section 3.2 shows that after the adding  formally $p$-th roots $(\tilde{x}_i)_{i=1,\dots, 2n}$ of central elements $({\xh}_i^p)_{i=1,\dots, 2n}$,
we identify (the pullback of) the $n$-th Weyl algebra
 with the algebra of matrices of size $(p^n\times p^n)$. Operator $P$ gives a matrix $M_{P,p}$ with coefficients being polynomials of a bounded degree independent on $p$ in generators  $(\tilde{x}_i)_{i=1,\dots, 2n}$. 
The support at prime $p$ is the hypersurface given by the equation 
  $$\det(M_{P,p})=0$$
 and hence has degree $const \cdot p^n$ in $(\tilde{x}_i)_{i=1,\dots, 2n}$. This polynomial has zero derivative with respect to each variable, and therefore is in fact a polynomial in $(\tilde{x}_i^p=\xh_i^p)_{i=1,\dots, 2n}$,
 of degree $const \cdot p^{n-1}$. 

The computer experiments leave no doubt that for ``typical'' $P$ this polynomial has indeed such a large degree, and
also is indecomposable. Hence we can not make a reasonable algebraic limit as $p\to \infty$.

Nevertheless, for certain operators $P$ the resulting polynomial is $p^{n-1}$-st power of a polynomial
 whose degree is uniformly bounded in $p$.
 This happens, for example, for the Hamiltonian of the periodic Toda lattice
$$\sum_{i\in \Z/n\Z}\left( \left(\frac{y_i\partial}{\partial y_i}\right)^2+\frac{y_{i+1}}{y_i}\right)$$
which after the transcendental change of variables $x_i=\log y_i$ has a more familiar form
  $$\sum_{i\in \Z/n\Z}\left( \left(\frac{\partial}{\partial x_i}\right)^2+\exp(x_{i+1}-x_i)\right)\,.$$

In general, we conjecture that such a situation is related to integrable systems:
\begin{conj} {\rm For a differential operator $P\in \D(X)$, the  support at prime $p$ of $\D(X)/\D(X)\cdot (P+\lambda\cdot 1)$
 is divisible by $p^{n-1}$ for generic constant $\lambda$ if and only if $P$ belongs to a quantum integrable system, i.e., 
 $P$ belongs to a finitely generated commutative $\k$-subalgebra of $\D(X)$ of Krull dimension $n=\dim X$.}
\end{conj}

The modification of $P$ by a generic additive constant is necessary in order to exclude certain parasitic examples, e.g.
$P=P_1 P_2$ where $P_1$ and $P_2$ belong to two different integrable systems.

\subsection{Donagi-Markman construction}

In \cite{Donagi}  a  construction of integrable systems was proposed
starting from any smooth projective variety $X$. One considers the scheme $B$ parametrizing smooth projective Lagrangian submanifolds $L\subset T^*X$ of a bounded degree. Scheme $B$ is smooth of dimension equal to $\dim \Gamma(L,T^*_L)$.
 Now consider the bundle $M$ over $B$ with the fiber over $[L]$ equal to the Albanese variety of $B$ (i.e., the abelian variety dual to $\op{Pic}_0(B)$). Then $M$ carries a natural symplectic structure,
 and fibers over $B$ are Lagrangian abelian varieties. Hence, we obtain an algebraic completely integrable system. In the special case $\dim X=1$
 one get a Zariski open part in the Hitchin integrable systems for the group $GL(N)$ for  $N:=\op{deg}(L\to X)$.
 
Smooth logarithmic families of non-compact Lagrangian submanifolds provide a natural generalization of Donagi-Markman construction.
The Albanese variety one can replace by an appropriate logarithmic version.

In \cite{Kontsevich2}, section 2, one can find a new notion of an algebraic integrable system which includes  simultaneously both the classical
and the quantized cases. This notion makes sense also over finite and local fields. The program in \cite{Kontsevich2}
was proposed relating integrable systems
 and Langlands correspondence in the functional field case. Roughly speaking, the spectrum of the maximal  commutative subalgebra in the quantized system should
 parametrize in certain sense holonomic $\D_X$-modules. Although the subject is not yet fully developed, it looks  certain that
 the integrable systems which should appear in the context of this program are the same as generalized Donagi-Markman systems.

Finally, we should mention that logarithmic families of planar curves appear naturally in matrix models. For example, one of the standard families in matrix theory (see \cite{DijkgraafVafa})
is the family of hyperelliptic curves
$$y^2=(W'(x))^2+ f(x),\,\,\,f(x)= \sum_{i=0}^{n-1} c_i x^i$$ where $W(x)$ is a fixed polynomial of degree $(n+1)$ and $(c_i)_{i\le {n-1}}$ are parameters of the  curve.
 The full solution  of the matrix model is based on the associated family of holonomic $\D_{\A^1_\k}$-modules.
 Also,  
 Seiberg-Witten curves for $N=2$ pure $SU(2)$ gauge theory, given by the equation
   $$x+1/x=y^2+u   $$
   where $u$ is a parameter, form a logarithmic family of curves in $T^*{\bf G}_m$ 
 endowed with coordinates $x,y\,\,\,\,(x\ne 0)$ and the symplectic form $dx \wedge dy/x$.

\section{Trigonometric version}
\subsection{Quantum tori}
The multiplicative (trigonometric) analog of the Weyl algebra $A_n$ 
is the algebra of functions on quantum torus:
$$K\langle \xh_1^{\pm 1},\dots,\xh_{2n}^{\pm 1}\rangle/\,(\xh_i\cdot\xh_j=q^{\omega_{ij}}\xh_j\cdot\xh_i)$$
where $\omega=(\omega_{ij})$ is the standard skew-symmetric matrix. The ground field is 
 $$K:=\k(q)\,,$$
 the field of rational functions in variable $q$, with coefficients in an algebraically closed field $\k=\overline{\k},\,\,\mathsf{char}(\k)=0$. There  exists a theory of holonomic modules over quantum tori, see \cite{Sabbah}.
 Examples of such modules are non-trivial cyclic modules in the case $n=1$, and also bimodules corresponding to automorphisms of skew-fields of fractions.
  An analog of a prime here is a primitive root of $1$.
   Namely, for any $N\ge 1$ and $q_N\in \k$ a primitive $N$-th root of 1,
    the algebra 
    $$\k\langle \xh_1^{\pm 1},\dots,\xh_{2n}^{\pm 1}\rangle/\,(\xh_i\cdot\xh_j=q_N^{\omega_{ij}}\xh_j\cdot\xh_i)$$
     has a large center equal to
      $$\k[ (\xh_1^N)^{\pm 1},\dots,(\xh_{2n}^N)^{\pm 1}]$$
       and it is an Azumaya algebra over its center, a twisted form
        of the algebra of matrices of size $N^n\times N^n$.
      
      In the case $n=1$  any element 
      $$P=\sum_{i,j\in\Z,\,|i|+|j|\le d} a_{ij}\xh_1^i \xh_2^j $$
      gives  a sequence of $N\times N$-matrices for $q_N\in \k$ a primitive $N$-th root of 1
      and $N$ large enough. Their determinants can be calculated similarly to section 3.3.
      The resulting algebra of ``periods'' is close to the algebra of expressions which appear e.g. as quantum Chern-Simons invariants of 3-dimensional manifolds.

 Here is an analog for quantum tori of the main conjecture from \cite{KaKo}.
 \begin{conj} 
{\rm There exists a homomorphism from the group $BirSympl_{n,\k}$ of birational symplectomorphisms the algebraic torus ${\mathbb G}_{m,\k}^{2n}$   endowed with the standard symplectic form $\sum_{i,j\le 2n} \omega_{ij} (x_i^{-1}dx_i)\wedge (x_j^{-1}dx_j)$, to the group of 
outer automorphisms of the skew field of fractions of the quantum torus. Also, the semiclassical
 limit as $q\to 1$ exists and gives the identity map from 
 the group of birational symplectomorphisms the group of birational symplectomorphisms the algebraic torus to itself.}
\end{conj}

    Let us assume the conjecture. Then, taking the reduction at a root of one $q=q_N$,  we  obtain an outer birational automorphism
 of the algebra with a large center, hence an automorphism of the center.
In this way one should obtain some mysterious ``quantum Frobenius'' endomorphism 
$$\op{Fr}_N: BirSympl_{n,\k}\to BirSympl_{n,\k}\,.$$
 There is a class of birational transformations analogous  to tame transformations (see \cite{KaKo}) in the additive case. It is generated by $Sp(2n,\Z)$,  by multiplicative translations
        $$M_{({\lambda_i})_{i\le 2n}}:(x_i)_{1\le i\le 2n}\mapsto (\lambda_i x_i)_{1\le i\le 2n}$$ for some constants
        $\lambda_i\in \k^*, 1\le i\le 2n$, and by the following non-trivial
         automorphism:
         $$T:x_1\mapsto x_1 (1-x_{n+1}),\,\,x_i\mapsto x_i\mbox{ for }i\ge 2\,.$$
          The miracle is that for the corresponding quantum automorphism
         $$\xh_1\mapsto \xh_1 (1-\xh_{n+1}),\,\,\xh_i\mapsto \xh_i\mbox{ for }i\ge 2$$
          the induced transformation on the center for $q=q_N$, is 
          given by the same formula:
       $$\xh_1^N\mapsto \xh_1^N (1-\xh_{n+1}^N),\,\,\xh_i^N\mapsto \xh_i^N\mbox{ for }i\ge 2\,.$$  

The quantum Frobenius $\op{Fr}_N$ acts identically on  $Sp(2n,\Z)$ and on transformation $T$, but on multiplicative translations
  it acts non-trivially
$$\op{Fr}_N :\,M_{({\lambda_i})_{i\le 2n}}\mapsto M_{({\lambda_i^N})_{i\le 2n}}\,.$$
We see that graphs of elements of the group $\Gamma_n$ generated by $Sp(2n,\Z)$ and $T$ are multiplicative analogs of isolated points
 (see section 4.4), i.e., could be considered as the ``most canonical'' Lagrangian subvarieties (and the corresponding holonomic modules).
The group $\Gamma_n$ deserves further study, it contains both all 
arithmetic groups and mapping class groups for large $n$, as follows e.g. from the work of Goncharov and Fock
 on cluster transformations and generalizations of Penner coordinates, see \cite{FockGoncharov}. Finally, the above discussion has a generalization
to elliptic algebras.

\vspace{5mm}

IHES, 35 route de Chartres, Bures-sur-Yvette 91440, France

{maxim@ihes.fr}
 
\end{document}